\documentstyle[11pt]{article}

\newtheorem{de}[th]{Definition}

\newtheorem{lem}[th]{Lemma}
\newtheorem{co}[th]{Corollary}
\newtheorem{re}[th]{Remark}

{
\newtheorem{th}{Theorem}[section]

\begin{document}

\large
\title {\bf METRIC AND GENERALIZED PROJECTION OPERATORS IN BANACH
SPACES:\\ PROPERTIES  AND APPLICATIONS}
\author{\bf Ya. I. Alber\\
Department of Mathematics\\
Technioin-Israel Institute of Technology\\
Haifa 32000, Israel}
\date{}
\maketitle

\section{Introduction}
\setcounter{equation}{0}
Metric projection operators in Hilbert and Banach spaces
are widely used in different areas of mathematics such as
functional and numerical analysis, theory of optimization
and approximation and for the problems of optimal control
and operations research, nonlinear and stochastic programming and games
theory.

Metric projection operators  can be defined in similar way
in Hilbert and Banach spaces. At the same time, they differ
signifitiantly in their properties.

Metric projection operator in Hilbert space is a monotone
(accretive) and nonexpansive operator. It provides an absolutely best
approximation for arbitrary elements from Hilbert space by the
elements of convex closed sets .
This leads to a variety of applications of this operator
for investigating theoretical questions in analysis and for
approximation
methods.  (For details see \cite{breg,pet,zd,li,aus,ct,al1,al2}).

Metric projection operators in Banach space do not have properties
mentioned above and their applications are not straightforward.
(See \cite{r,aln3}).

On the other hand, questions of smoothness and differentiability
of metric projection operators in Banach spaces were actively
investigated \cite{bj,li,aln2,al2}. New results in this field are
immiadetely used
in various applications. For example, recently established
in \cite{aln2,al2} properties of uniform continuity of these operators
were used in \cite{al3,al4} to prove  stability  of  the  penalty  and
quasisolution methods.

Two of  the most important applications of the method of metric
projection
operators are as follows:
\begin{itemize}
\item  solve a variational inequality by the iterative-projection method,
\item  find common point of convex sets by the iterative-projection method.
\end{itemize}
(See Section 2 and Section 3 below).

In Banach space these problems can not be solved in the framework of metric
projection operators. Therefore,
in the present paper we introduce  new generalized projection operators
in Banach space as a natural  generalization of  metric projection
operators in Hilbert space.  To demonstrate our approach, we apply
these operators for solving two problems mentioned above in Banach
space.

In Section 2 and Section 3 we introduce notations and recall some results from
t
theory of variational inequalities and theory of approximation. Then
in Section 4 and Section 5 we describe the properties of metric projection
operators  $P_\Omega$ in
Hilbert and Banach spaces and also formulate equivalence theorems
between variational inequalities and direct projection equations  with these
operators.
In Section 6 we discuss the existence of strongly
unique
 best approximations  based on Clarkson's and parallelogram  inequalities.

In Section 7 we introduce generalized projection operator $\Pi_\Omega$ which
acts from Banach space $B$ on  the convex closed set $\Omega$ in the
same
space $B$. Then we state its properties and give the convergence theorem
for the
method of successive generalized projections  used to find a common
point of
convex sets.

In Section 8 new generalized projection operator $\pi_\Omega$
acting from  conjugate Banach space $B^{*}$    on convex closed set
 $\Omega$ in  the space $B$ and its properties are examined.
Then a theorem of equivalence of
the solutions of variational inequality and operator equation with
operator $\pi_\Omega$ is presented.  It constitutes the basis for
construction of iterative-projection methods for nonlinear problems in
Banach space
(including smooth and nonsmooth optimization problems)

Note that the main properties of metric and generalized projection operators
in Banach spaces have been obtained by using principally new technique
including Banach space geometry, parallelogram inequalities, nonstandard
Lyapunov functionals and estimates of moduli of monotonicity and continuity
for duality mappings.

\section{Variational Inequalities. Problems of Optimization}
\setcounter{equation}{0}
In this Section we recall some of the results from the
theory of variational inequalities and formulate a problem on the
equivalence
between solutions of the variational inequalities and corresponding
operator
equations.  These equations are solved by the iterative-projection
methods. This
yields an aproximation of solutions of the initial variational
inequalities.

Let $B$ be a real (reflexive) uniformly convex and uniformly smooth
Banach
space, $B^{*}$ its conjugate (dual) space, $ ||\cdot ||, ||\cdot
||_{B^{*}},
||\cdot ||_H$ norms in the Banach spaces $ B, B^{*}$ and in Hilbert
space $ H$.
 As usually we denote $<\varphi,x>$ a dual  product in $B$ and.
 This determines pairing between $\varphi  \in B^{*}$ and $x \in  B$
\cite{br} ($(y,x)$ is an inner product in $H$).
Let $\Omega $ be a  nonempty convex closed set in B.

\begin {de}  \label{gpo}
The operator $P_{\Omega}: B \rightarrow {\Omega} \subset B$
is called metric projection operator if it yields
the correspondence between an arbitrary point  $x \in  B$ and nearest
point
 $\bar x \in \Omega $ according to minimization problem
\begin{equation} \label{f1}
P_{\Omega}x = \bar x; \;\;\; \bar  x: ||x - \bar x|| =
\inf_{\xi \in \Omega} ||x - \xi||.
\end{equation}
 \end {de}

Under our conditions operator $P_{\Omega}$ is defined at any point $x
\in  B$ and
it is single-valued, i.e. there exists for each $x \in  B$ a unique
$projection $
$ \bar x$  called best approximation \cite{h}.

Let $A$ be an arbitrary operator acting from Hilbert space  H
to H ,
$\alpha$ an arbitrary fixed positive number  and ($\varphi,x$)
an inner product in H.  Let also $f \in H$.
It is well  known that (see, for instance,  \cite{cea})
\begin {th} \label{th1}
 The point $x \in \Omega  \subset  H$ is a solution of variational
inequality
\begin {equation} \label{f2}
(Ax-f, \xi - x) \ge 0, {\qquad }        \forall \xi \in\Omega
\end{equation}
 if and only if x is a solution of the following operator equation in H
\begin{equation} \label {f3}
x=P_{\Omega}(x - \alpha (Ax-f)).
\end{equation}
\end{th}

This is an important statement, because it provides a basis for
constructing
approximate (iterative) methods  in Hilbert spaces. The simplest  method
of this type can be described as follows
\begin{equation} \label{f4}
x_{n + 1} = P_{\Omega}(x_{n}-\alpha_n (Ax_{n} - f)), \qquad n=0,1,2...
\end{equation}
Under suitable choice of relaxation parameters $\alpha_n$ , one can
prove
that iterative process (\ref {f4}) converges strongly to the solution of
the
variational inequality (\ref {f2}).
It can be done for operator $A$ which have different structures and
different types of smoothness  \cite{al1,cea,aus,dunn}.
Moreover, one can establish both stability and nonasymptotic estimates
of convergence rate \cite{al1}.

     We want to  emphasize that  the problem of solving operator
equation
$Ax=f$ and the problem of minimization of the functional $u(x)$ on
$\Omega$ are
 realized as variational inequalities (\ref {f2}) for $\Omega =  H$ and
$$Ax=\partial u(x),{\quad } x \in\Omega, {\quad } f = 0 $$
respectively. Here $\partial u(x)$  is gradient  or subgradient of the
functional
 $u(x)$.

Now we consider more general and more complicated case
of the variational inequality
\begin {equation} \label{f5}
<Ax - f, \xi - x> \ge 0, {\qquad }       \forall \xi \in\Omega
\end{equation}
in Banach space $B$ with operator $A$ acting from  $B$ to $B^{\ast}$
\cite{aln3}.
There is a natural problem to formulate and to prove  an analogue of
Theorem \ref{th1} in Banach space, and then to use it as a basis to
construct
iterative-projection methods similar to (\ref{f4}).

It is quite obvious that  the Banach space analogue of the equation
(\ref{f3}) has the
following form
\begin {equation} \label{f6}
x=\Gamma_{\Omega}(Fx-\alpha(Ax-f))
\end {equation}
with operator F  acting from  B to $B^{\ast}$. The equation (\ref{f6})
is unusual  because operator ${\Gamma}_\Omega$  "projects" elements from
the dual space $B^{\ast}$ on the set $\Omega \subset B$.
Here one can not use metric projection
operator $P_{\Omega}$ for this purpose because it acts from  B
to B.  It turned out that a natural generalization of metric
projection operator in  Hilbert space leads to a new operator which
we call generalized projection operator:
 $$\pi_{\Omega}  : B^{\ast} \rightarrow  \Omega \subset B .$$
This automaticly yields the following form of the equation (\ref{f6})
 $$x= \pi_{\Omega } (Jx-\alpha(Ax-f))$$
where $J: B\rightarrow B^{\ast}$ is a normalized duality mapping in $B$
\cite{br}.
The operator $J$ is one of the most significant operators  in nonlinear
functional
analysis. In particular, it is used in the theory of optimization and in
the theory of
monotone and accretive operators in Banach spaces.  It is determined by
the
expression
 $$ <Jx,x> = ||Jx||_{B^{*}}||x|| = ||x||^2.$$
 Note also that a duality mapping exists in each Banach space.
 In what follows we recall from \cite{al3} some of the examples of this
mapping
in spaces
$l^p, L^p $ , $ W^p_m$ , $\infty >p>1 $ :
\begin{itemize}
\item  (i) $ l^p:  Jx = ||x||^{2-p}_{l^p} y\in l^q, {\quad } x =
\lbrace  x_1, x_2,...\rbrace  $,
 $ y =   \lbrace  x_1 ||x_1||^{p-2}, x_2 ||x_2||^{p-2},...\rbrace, \\
 p^{-1} + q^{-1} = 1 $
\item  (ii) $ L^p:  Jx = ||x||^{2-p}_{L^p} |x|^{p-2}x \in L^q $
 \item  (iii) $ W^p_m:  Jx = ||x||^{2-p}_{W^p_m} \sum (-1)^ {|\alpha
|}D^\alpha
   (|D^\alpha x |^{p-2} D^\alpha x ) \in W^q_{-m} $
 \end{itemize}
 Note that in Hilbert space  $J$ is  an identity operator.

 Now we  define the iterative method similar (\ref{f4}) as follows
\begin {equation} \label{k6}
x_{n + 1} = \pi_{\Omega}(Jx_{n}-\alpha_n (Ax_{n} - f)), {\qquad
}n=0,1,2...
\end{equation}
We will give full mathematical foundation for this method including
three basic
aspects: convergence, stability and estimates of convergence rate, in
forthcoming
paper.

\section{ Problems of Approximation. Common Points of Convex Sets}
\setcounter{equation}{0}
Second important problem which is investigated in this paper  using
projection operators can be formulated as follows:  find  common point
of an
ordered collection of convex and closed  (i.e. Chebyshev) sets  $\lbrace
\Omega_{1},
\Omega_{2 },...,\Omega_{ m} \rbrace $ in uniformly  convex Banach  space
$B$.
Here we assume that  sets $\lbrace  \Omega_{i},  i=1,2,...,m \rbrace$
have nonempty intersection  $\Omega _{*}= \bigcap _{i=1}^{m}
\Omega_{i}$.
Let us define a composition
\begin {equation} \label{f9}
P =  P_{1} \circ  P_{2} \;\circ \cdotp \cdotp \cdotp \circ \;
 P_{m},  \qquad     P_{i} = P_{\Omega_{i}}
\end{equation}
 and introduce  method of successive projections according to a formula
\begin {equation} \label{f10}
x_{n+1} = P^{n+1}x_0,  {\qquad }     n=0,1,2,...,  {\qquad } x_o \in B.
\end{equation}

Convergence of the iterative process  (\ref {f9})  and  (\ref {f10}) as
well as of similar
processes to the  point  $ x_{*} \in \Omega _{*} $ was proved before
only in Hilbert
space .  (See \cite{vn,breg,bc,ct,d,gpr}).  In the formulae (\ref {f9})
and
 (\ref {f10})
which describe method of successive projections one can use metric
projection
operators $P_{\Omega}$ in Banach space. However, up to this date
no proof was suggested for the convergence of (\ref {f9}) and (\ref
{f10})
in Banach space. The reason is that in Hilbert space H the metric
projection
operator satisfies the following  significant inequality
\begin {equation} \label{f3.3}
||P_\Omega{ x} - \xi||_H \leq||x-\xi||_H , {\qquad }        \forall \xi
\in\Omega
\end{equation}
which can be obtained  from the property of nonexpansiveness of  this
operator
in H
 \begin {equation} \label{f3.4}
 {||P_\Omega{x} - P_\Omega{y}||_H} \leq {||x- y||_H}.
\end{equation}
 It satisfies also a much stronger property (see Section 4)
\begin {equation} \label{f13}
{||P_\Omega{ x} - x||^2}_H
\leq {||x-\xi||^2}_H - {||P_\Omega{ x} - \xi||^2}_H ,  {\qquad }
\forall \xi \in\Omega .
\end{equation}
But in Banach spaces these properties do not hold in general
\cite{sm,r,gr}.

For example, in \cite{sm}  it is shown that in uniformly convex  Banach
space
with modulus of convexity   $\delta (\epsilon) $ of order  $\epsilon ^q
$,
$q\ge 2 $ \cite{dis,lz,ml},  the inequality
\begin {equation} \label{f14}
{||P_\Omega{ x} - x||}^q \leq {||x-\xi||}^q - \lambda  ||P_\Omega{ x} -
\xi||^q , {\qquad }
        \forall \xi \in\Omega .
 \end{equation}
holds. Coefficient $\lambda < 1$ in (\ref{f14}) and it depends on $q$.
Namely, in \cite{sm} it is defined in Banach spaces of the type
$ L^p$ and $ W^p_m,\;1 < p < \infty ,$ as follows (cf. (\ref {f161}) and (\ref
{
$$ 1 < p \leq 2,\;\; q=2 , \; \;  \lambda = (p-1)/8;$$
and
$$ 2 < p < \infty,\;\; q=p , \;  \;  \lambda = 1/p2^p.$$

Inequality (\ref{f14}) yields
\begin {equation}  \label{f15}
{||P_\Omega{ x} - \xi||}^q \leq {\lambda}^ {-1}{||x-\xi||}^q -
{\lambda}^ {-1}{||P_\Omega{ x} - x||}^q,
{\qquad }  \forall \xi \in\Omega.
\end{equation}
This does not guarantee the nonexpansiveness of   the metric projection
operator in
Banach space even for $\xi = y\in\Omega$  while in (\ref {f3.4}) $\xi =
y $ is
an arbitrary element of the space H. But without this property, one can not
to study the method (\ref {f9}), (\ref {f10}).

Now  we consider more general case.
In \cite{sm} a strongly unique best approximation  was defined as
follows.
\begin{de} \label{def1}
$\bar x$ is called a strongly unique best approximation
in $\Omega$ for the element  $x  \in  B$ if there exists a constant
$\lambda$
and a strictly increasing function $\phi (t) : R^{+}  \rightarrow R^{+}$
such that
 $\phi(0)=0$ and
\begin {equation}  \label {f151}
{\phi(||P_\Omega{ x} - x||)} \leq {\phi (||x-\xi||)} - \lambda \phi
(||P_\Omega{ x} - \xi||) ,
{\qquad }  \forall \xi \in\Omega.
\end{equation}
\end{de}

The projection $\bar x = P_\Omega  x$ in Hilbert space, and the
projection
$\bar x$ in Banach space under the conditions of  \cite{sm}, are
strongly unique best
 approximations in $\Omega $.

We  call the projection $\bar x $ {\it absolutely best approximation \/}
of $x\in B$
with respect  to  function $\phi (t) $ if  $\lambda  = 1$ in (\ref
{f151}).
In this case the inequality (\ref {f151}) can be represented in the
equvalent form
 $$\phi (||P_\Omega{ x} - \xi||) \leq {\phi (||x-\xi||)} - \phi
(||P_\Omega{ x} - x||),
 {\qquad }  \forall \xi \in\Omega. $$
 It is clear from (\ref {f13}) that projection $\bar x$ in Hilbert space
is absolutely
best approximation with respect to  function  $\phi (t) = t^2 $ (or with
respect to
functional $\phi (\xi) =  {||x-\xi||^2}_H $ with  $x$  fixed). But it is
not true in Banach spaces.

So, metric projection operator can not be used in (\ref {f9}) and
(\ref {f10}).
Instead  we introduce new generalized projection operator
$$\Pi_{\Omega}:  B \rightarrow {\Omega} \in B $$
so that the inequalities (\ref{f3.3}) and (\ref{f13}) hold for  some
Lyapunov
functional $\phi (\cdot)$. In Section 7 we will provide the convergence
theorem for the process
(\ref {f9}) and (\ref {f10}) which now has the form
\begin {equation} \label{f16}
x_{n+1} = \Pi^{n+1}x_0,  {\qquad }     n=0,1,2,..., {\qquad } x_o \in B
\end{equation}
and
\begin {equation}  \label{f17}
\Pi =  \Pi_{1} \circ  \Pi_{2} \;\circ \cdotp \cdotp \cdotp \circ \;
 \Pi_{m},  \qquad      \Pi_{i} = \Pi_{\Omega_{i}}.
 \end{equation}

\section{ Metric Projection Operator  $P_{\Omega}$ in Hilbert Space}
\setcounter{equation}{0}
All results described in \S 2 and \S 3 for two basic problems were
obtained
only in Hilbert space. This is due to the fact that many remarkable
properties
of the metric projection operators can not be extended from Hilbert
space to
Banach space. This is way we introduce in \S 7 and \S 8 new generalized
projection operators in Banach spaces which have all properties of
metric projection operators in Hilbert space .

Before that in \S 4 and \S 5 we compare complete lists of the properties
of the
metric  projection operators in Hilbert  and Banach spaces.

We denote $ \bar x = P_{\Omega}x $. Let $\xi \in {\Omega}$ and $
{\Omega}
\subset H$.
The following properties are  valid in Hilbert space:
\cite {aus,breg,cea,gr,h,pet,al3,n,z}.

{\bf 4.a.}  $P_{\Omega}$ is fixed at each point $\xi $, i.e.
 $P_{\Omega}\xi = \xi.$

{\bf 4.b.}  $P_{\Omega}$ is monotone (accretive) in $ H$, i.e.
$$(\bar x - \bar y, x-y) \ge 0 .$$

{\bf 4.c.} $ (x - \bar x, \bar x -\xi) \ge 0,  {\qquad }        \forall
\xi \in\Omega   .$
The point $\bar x$ is metric projection of $x$  on $\Omega \subset H$
iff the inequality 4.c is satisfied.
The property 4.c we call {\it basic variational principle} for $P_{\Omega}$
in $H$.

{\bf 4.d.} $ (x - \xi, \bar x -\xi) \ge 0,  {\qquad }        \forall \xi
\in\Omega  .$

{\bf 4.e. } $ (x - \bar x,  x -\xi) \ge 0,  {\qquad }        \forall \xi
\in\Omega .$\\
In fact, even stronger inequality
$$ (x - \bar x,  x -\xi) \ge  ||x - \bar x||^2_H ,  \;\; \forall \xi
\in\Omega$$
holds  (see (\ref {f52})).

{\bf 4.f.}   $P_{\Omega}$ is nonexpanseve in $ H$,  i.e.

 $${ ||\bar x - \bar y||_H} \leq {||x- y||_H} .$$

{\bf 4.g.}  $P_{\Omega}$ is $P$-strongly monotone in $ H$,  i.e.
 $$(\bar x - \bar y,  x-y) \ge  ||\bar x - \bar y||^2_H .$$

{\bf 4.h. } The operator $P_{\Omega}$ yields an absolutely best
approximation
of $x\in H$ with respect to the functional  $ {V_1 (x,\xi )=  ||x -
\xi||^2 _H}$

$${||\bar x  - \xi||}^2 _H \leq {||x-\xi||}^2 _H - {||x - \bar x||}^2_H
,
{\qquad }
\forall \xi \in\Omega .$$

{\bf 4.i.} Any $P_{\Omega}$ satisfies the inequality (see \cite {n})
$$((I - P_{\Omega})x - (I - P_{\Omega})y, P_{\Omega}x - P_{\Omega}y) \ge 0,
\;\;\forall x,y \in H.$$

\section{ Metric Projection Operator  $P_{\Omega}$ in  Banach Space}
\setcounter{equation}{0}
Here we show that some of the properties of the metric projection
operators in Hilbert
space are not satisfied in Banach space. At the same time, we describe
in detail the
properties  of uniform continuity of the metric projection operators in
Banach space.

We denote $ \bar x = P_{\Omega}x $. Let $\xi \in {\Omega}$ and $
{\Omega}\subset B$.
The following properties hold in Banach space (the sign "-" from
5.a - 5.i denotes an absence of corresponding property):

{\bf  5.a.} The operator $P_{\Omega}$ is fixed at each point $\xi $, i.e.
 $P_{\Omega}\xi = \xi .$

{\bf  5.b.} -

{\bf  5.c.}  $ < J(x - \bar x), \bar x -\xi > \ge 0,  {\qquad }
\forall \xi \in\Omega  $ (see \cite{li}).
The point $\bar x$ is metric projection of $x$  on $\Omega \subset B$
iff the inequality 5.c is satisfied.
The property 5.c we call {\it basic variational principle} for $P_{\Omega}$
in $B$.

{\bf  5.d.} -

{\bf  5.e. } $ < J(x - \bar x),  x -\xi > \ge 0,  {\qquad } \forall \xi
\in\Omega $.  \\
It follows from the next stronger statement .
\begin {th} \label{f51}
$\bar x \in \Omega $ is a projection of the point $x \in B$ on $\Omega$
if and only
if the inequality
\begin {equation} \label{f52}
  < J(x - \bar x),  x -\xi > \ge ||x - \bar x||^2 ,  {\qquad }  \forall
\xi \in\Omega
\end {equation}
is satisfied.
\end {th}

{\bf Proof} In fact, from (\ref {f52}) it follows immediately that
$$||x - \bar x|| \leq  {||x - \bar x||}^{ -1} { < J(x - \bar x),  x -\xi
> } \leq  ||x - \xi||,
\qquad   \forall \xi \in \Omega $$
 i.e.  $\bar x = P_{\Omega}x $.  Inversely, if   $ \bar x = P_{\Omega}x
$,  then by
virtue of  5.c  we have
\begin{equation}
 0 \leq {< J(x - \bar  x), \bar  x -\xi > }
 =  {< J(x - \bar  x), \bar  x - x > } +  {< J(x - \bar  x),  x -\xi > }
 =
\end{equation}
$$ = -||x - \bar  x||^2  +   < {J(x - \bar  x),  x -\xi > }$$
which yields (\ref {f52}).

{\bf  5.f.} Now we describe the property of uniform continuity of
operator $ P_{\Omega}x $
in Banach space $B$. Recall that in  Banach space  the metric projection
operator  is not nonexpansive in general case. But it is uniformly
continuous on
each bounded set according to the following theorem.
\begin {th} \label{f54}
 Let B be an uniformly convex and  uniformly  smooth Banach space.  If
$\delta _B (\epsilon) $  is a modulus of convexity of the space B,
$g_B (\epsilon)  = \delta  _B (\epsilon)  / \epsilon$  and  $g^{-1}_B
(\cdot)$
is an inverse function, then
\begin{equation} \label{f55}
||\bar x - \bar y|| \leq  C  g^{-1}_B  (2LC^2 g^{-1}_{B^{*}} {(2CL||x -
y||)}),
\end{equation}
where $1 < L  < 3.18$ is Figiel's constant (see \cite{fg}) and
$$ C = 2  {max} \lbrace 1, ||x - \bar y ||, ||y - \bar x||\rbrace .$$
\end {th}
\begin{re}  If  $ ||x - \bar y || \leq R$ and $||y - \bar x|| \leq R$,
then
 $(C = 2 {max} \lbrace 1, R \rbrace) $ is an absolute constant and
(\ref{f55}) is a
provides a quantitative description of the uniform continuity  of
operator $ P_{\Omega} $
in Banach space on each bounded set.
\end{re}
The estimate (\ref{f55}) which was established in \cite{aln2} is global
in nature.
Earlier, in \cite{bj}  Bjernestal obtained local estimate
\begin {equation} \label{f57}
 ||\bar x - \bar y|| \leq 2\delta^{-1}_B (6\rho _B {(2||x - y||)}),
\end{equation}
 where $\rho _B (\tau)$ is  a modulus of smoothness of the space B
\cite{dis,lz,ml} .

The estimate of (\ref{f57}) is better than  our estimate (\ref{f55}).
This is why in
 \cite{al2} we continued the investigation of uniform continuity of
metric
projection operator in Banach space.  It turns out that the following
global variant
of (\ref{f57}) can be obtained.
\begin {th} \label{f58}
 Let B be an uniformly convex and  uniformly  smooth Banach space. If
$\delta _B (\epsilon) $  is a modulus of convexity of the space $B$ and
 $\rho _B (\tau)$ is  a modulus of  its  smoothness, then
\begin {equation} \label{f59}
||\bar x -\bar y|| \leq C\delta^{-1}_B (\rho _B {(8CL||x - y||)}),
\end{equation}
where  constant L and  function C are defined in Theorem \ref{f54}.
\end {th}

\begin{re}  To accuracy constants  the estimates (\ref{f55}) and
(\ref{f59}) give
respectively
 $$||\bar x - \bar y|| \leq  g^{-1}_B  ( g^{-1}_{B^{*}} {(||x - y||)}),
$$
 $$||\bar x - \bar y|| \leq \delta^{-1}_B (\rho _B {(||x - y||)}). $$
\end{re}

{\bf  5.g.} -

{\bf  5.h.} -

{\bf 5.i.} Any $P_{\Omega}$ satisfies the inequality (see \cite {li})
$$<J(x - P_{\Omega}x) - J(y - P_{\Omega}y), P_{\Omega}x - P_{\Omega}y> \ge 0,
\;\;\forall x,y \in B.$$

Using  the properties of  metric projection operator $P_\Omega$
we obtained
Banach space analogue of Theorem \ref{th1}.
\begin {th} \label{f580}
 Let $ A $ be an arbitrary operator from Banach space $B $ to $B^{*}$,
$\alpha$ an
 arbitrary fixed positive number, $f\in B^{*}$.  Then  the point $x \in
\Omega  \subset B$
is a solution of variational inequality
 $$<Ax-f, \xi - x> \ge 0, {\qquad }        \forall \xi \in\Omega ,$$
if and only if x is a solution of the operator equation in B
\begin{equation} \label {f282}
x=P_{\Omega}(x - \alpha J^{*}(Ax-f))
\end{equation}
where  $ J^{*}: B^{\ast} \rightarrow B $ is normalized duality mapping
in $  B^{\ast}$.
\end{th}
{\bf Proof} Since $J$ is homogeneous and odd operator, and
$J J^{*} = I_{B^{*}}:B^{\ast} \rightarrow B^{\ast}, $
we can rewrite given variational inequality in form of
$$<(J(- \alpha J^{*}(Ax-f) + x) - x, x - \xi> \ge 0,
{\qquad }        \forall \xi \in\Omega. $$
The equation (\ref {f282}) is obtained immediately from
basic variational principle  for $P_{\Omega}$  in $B$ (the property 5.c).

Iterative process corresponding to (\ref{f282}) is the following
\begin{equation} \label {f2823}
x_{n + 1} = P_{\Omega}(x_{n}-\alpha_n J^{*} (Ax_{n} - f)),
\qquad n=0,1,2...
\end{equation}
However, there are not any approaches to investigation of (\ref{f2823}).
We conjecture that this iterative process is divergent.

In Section 7 and Section 8  we will construct the generalized projection
operators in
Banach spaces with the additional properties  5.b,  5.d,  5.g and  5.h.,
and we will provide  the equivalence theorems with such operators
and iterative methods for which one can establish convergence, stability,
and nonasymptotic estimates of convergence rate.

In conclusion of this Section we note that Theorem \ref{f580} can be proved
for more general variational inequalities considered in \cite{jy}.

\section{Parallelogram Inequalities and Strongly Unique Best
Approximations}
\setcounter{equation}{0}
In this section we discuss the existence of strongly unique best
approximation in
the spaces $l^p, L^p $  and $ W^p_m$  where  $\infty > p >1. $

In \cite{aln8} (see also \cite{aln9}) we established the following upper
parallelogram inequality
 $$2||x||^2 + 2||y||^2 - ||x + y||^2 \leq 4 ||x - y||^2 +
C_1 \rho _B (||x - y||) ,$$
$$ C_1 = C_1(||x||,||y||) = 2 \max \lbrace  L, (||x|| +||y||)/2 \rbrace $$
for the arbitrary points $x$ and $y$ from uniformly smooth Banach space
$B$.
We also obtained lower  parallelogram inequality  (see Section 7)
\begin {equation} \label{f167}
2||x||^2 + 2||y||^2 - ||x + y||^2 \ge L^{-1}\delta_B (||x -
y||/C_2,
\end{equation}
$$ C_2 = C_2(||x||,||y||) = 2 \max \lbrace  1, \sqrt{(||x||^2
+||y||^2)/2}\rbrace $$
for the arbitrary points   $x$ and $y$ from uniformly convex Banach
space $B$.
Analogous parallelogram inequalities for the $||x||^q$ of other orders
$q$
were obtained in \cite{not}.

If  $||x|| \leq R $ and $||y||\leq R $ then
$$ C_1(||x||,||y||) = \bar C_1 =
2  \max \lbrace  L, R\rbrace $$
and
$$C_2(||x||,||y||) = \bar C_2 =
2 \max \lbrace  1, R\rbrace. $$
The inequality (\ref {f167}) gives
the  uniform convexity of  the functional $||x||^2$ in the space
$B$ with modulus of convexity
$\delta (||x - y||) = (2L)^{-1}\delta_B (||x - y||/C_2).$
This can be established from the following lemma \cite{vnc}.
\begin {lem} \label{f77}
If a convex functional ${\phi(x)}$ difined on convex closed set $ \Omega $
satisfies the inequality
$${\phi({\frac 12}x + {\frac 12} y)} \leq {\frac 12} \phi (x) +
{\frac 12} \phi (y) -  \kappa (||x - y||),$$
where $\kappa (t) \ge 0,  \kappa (t_0) > 0 $ for some $t_0 > 0$,
then ${\phi(x)}$ is uniform convex of functional with modulus of convexity
$\delta (t) = 2 \kappa (t)$, and
$${\phi(x)} \ge  \phi (y) + <l(y),x - y> + 2 \kappa (||x - y||) ,$$
\begin {equation} \label{f109}
 <l(x) - l(y),x - y> \ge 4 \kappa (||x - y||)
\end{equation}
for all $l(x) \in \partial \phi(x)$. Here $ \partial \phi(x)$ is the
set of all support functionals (the set of all subgradients) of $\phi(x)$ at
the point $x \in \Omega.$
\end{lem}

{}From this Lemma and (\ref {f167}) it follows that
\begin {equation} \label{f168}
||x||^2  \leq ||y||^2 +  2<Jx, x-y>  -  (2L)^{-1}\delta_B (||x - y||/C_2).
\end{equation}
and
\begin {equation} \label{f268}
<Jx - Jy, x-y> \ge  (2L)^{-1}\delta_B (||x - y||/C_2).
\end{equation}
Let  $\Omega \subset B,\; \xi \in \Omega,\; \bar x = P_\Omega x $.  We
replace in (\ref {f168})  $x$ by $(x - \bar x)$ and $y$ by $(x - \xi)$
and
obtain
$$||x - \bar x||^2  \leq ||x - \xi||^2 - 2 <J(x - \bar x), \bar x - \xi>
- (2L)^{-1}\delta_B (||\bar x - \xi||/C_2).$$
The property  5.c  then yields the following general formula
\begin {equation} \label{f169}
 ||x - \bar x||^2  \leq ||x - \xi||^2  - \lambda \delta_B (||\bar x -
\xi||/C_2), \;\; \lambda = (2L)^{-1} .
\end{equation}
It is obvious  that  if $\delta_B (\epsilon)$ is proportional to
$\epsilon ^2$ or can be estimated by $\epsilon ^2)$
(this occures  in Hilbert spaces or in the spaces
of type $L^p$ for  $1 < p \leq 2 $ \cite{aln1},) then
projection $\bar x$ is a strongly
unique best approximation with $\phi (t) = t^2$, at least, on each bounded
set (See Def. \ref{def1}).  However, constant $\lambda$
in  inequality (\ref {f169}) is not exact in these cases.

Using the inequality (2.3) from \cite{lhh}, we immediately
obtain from Lemma \ref{f77} that in spaces $L^p$,  $1 < p \leq 2 $  the
estimate
\begin {equation} \label{f161}
||x - \bar x||^2  \leq ||x - \xi||^2  -  (p-1)||\bar x - \xi||^2  .
\end{equation}
is valid. This coincides with the result of \cite{lhh} (Theorem 4.1) which
has been obtained there in a different way.

In \cite{Yu} it was shown that in spaces   $W^p _m$ (consequently, in
$L^p $
and  $l^p$),  $1 < p \leq 2 ,$ the following inequality holds
\begin {equation} \label{g69}
 ||x + y||^2  \leq   2||x||^2 + 2||y||^2  - (p-1)||x - y||^2.
\end{equation}
Then Lemma \ref{f77}, the property 5.c and (\ref {g69}) give
$$||x - \bar x||^2  \leq ||x - \xi||^2  -  (p-1)||\bar x - \xi||^2 /2 . $$
This improves the estimate of a strongly unique best
approximation in the space $W^p_m$, $\infty > p > 2.$
(Compare with the corresponding inequality in \cite{sm}).

Furthermore, the strongly unique best approximation of the projection
$\bar x $  in spaces $B$ of the type  $l^p, L^p $  and $ W^p_m$  where
$\infty > p > 2, $
can be established from Lemma \ref{f77} and Clarkson's  inequality
$$ ||x+y||^p  \leq    2^{p-1} ||x||^p  +  2^{p-1}  ||y||^p  - ||x-y||^p.$$
This  means that functional $||x||^p$ is uniformly convex.
Therefore, we can write
\begin{equation} \label{f162}
 ||y||^p \ge ||x||^p +p <J^{\mu}x,y-x> + 2^{-p+1} ||x-y||^p
\end{equation}
where  $ J^{\mu}$  is a duality  mapping  with the gauge function $\mu (t) =
t^{p-1}$ (see  \S 7).
Now we substitute in  (\ref{f162}) $(x -\bar x)$ and $(x - \xi)$ for $x$
and $y$,
respectively.  By virtue  of
 $$<  J^{\mu} (x - \bar x), \bar x - \xi > \ge 0 $$
and
$$||x - \bar x||^p \leq  ||x-\xi||^p - p <J^{\mu}(x - \bar x),\bar x - \xi>
- 2^{-p+1}
||\bar x-\xi||^p $$
we have
\begin{equation} \label{f177}
||x - \bar x||^p \leq  ||x-\xi||^p  -    2^{-p+1} ||\bar x-\xi||^p .
\end{equation}
This improves the estimate of a strongly unique best
approximation even in the space $L^p$, $\infty > p > 2.$
(Compare with the corresponding inequality in \cite{sm,ls}).
It was obtained without any additional conditions on the set $\Omega$
and space $B$.

Besides, the important generalization of (\ref {f177}) is valid.
\begin {th}
Let $B$ be a space either $l^p$ or $ L^p $ or $ W^p_m$  where
$\infty > p \ge 2. $ Let $\Omega$ be a closed convex set in $B$. Then for
every point $x \in B$ there exists a unique point $\bar x = P_\Omega x $
such that
$$||x - \bar x||^s \leq  ||x-\xi||^s  -    2^{-s+1} ||\bar x-\xi||^s ,
{\quad} \forall \xi \in\Omega, {\quad} s \ge p \ge 2. $$
\end {th}

{\bf Proof} follows from the inequality (see \cite {vnc})
 $$||x+y||^s  \leq    2^{s-1} ||x||^s  +  2^{s-1}  ||y||^s  - ||x-y||^s,
  {\quad} s \ge p \ge 2. $$

\section{Generalized Projection Operator $\Pi_\Omega$ in Banach Space}
\setcounter{equation}{0}

Here we introduce generalized projection operator $\Pi_\Omega$  and
describe
its properties in Banach spaces. Then we formulate theorem about
convergence
of the method of succesive projections given in \S 3.  This method
yields
an approximation of the  common  point of convex closed sets. (See
second problem in \S 3).

The  formula (\ref{f1}) in the Definition \ref{gpo} of the metric
projection operator  is equivalent to the minimization problem
\begin{equation} \label{f61}
P_{\Omega}x = \bar x; \;\;\; \bar  x: ||x - \bar x||^2 =
\inf_{\xi \in \Omega} ||x - \xi||^2 .
\end{equation}

Note that $V_1(x, \xi) = ||x - \xi||^2$ can be considered not only as
square
of distance between points $x$ and $\xi$ but also as Lyapunov functional
with respect to $\xi$ with fixed $x$. Therefore,  we can rewrite (\ref
{f61}) in the
form
 $$P_{\Omega}x = \bar x; \;\;\; \bar  x:  V_1(x ,\bar x) =
\inf _{\xi \in \Omega} V_1(x ,\xi), {\qquad } $$
In Hilbert space
 $$V_1(x ,\xi) = ||x||^2 _H-2 (x, \xi) +||\xi||^2_H  .$$
In the papers \cite{aln8,aln1} we introduced  Lyapunov functional
\begin{equation} \label{k64}
V_2(Jx ,\xi) = ||Jx||^2_{B^{*}} -2 <Jx, \xi> +||\xi||^2 .
\end{equation}

 It is a nonstandard functional because it is defined on both the
elements $\xi$ from
 the primary space $B$ and on the elements $(Jx)$ from the dual space
$B^{*}$
(see also (\ref {f72})).
\begin {lem} \label{s1}
The functional $V_2(Jx ,\xi)$ has the following properties:\\

1.  $V_2$ is continuous.\\

2.  $V_2$ is differentiable with respect to $\varphi $ and $\xi$ .\\

3.  $ grad_{\varphi} V_2(Jx, \xi) = 2(x-\xi) .$\\

4.  $ grad_{\xi} V_2(Jx, \xi) = 2(J\xi-Jx) .$\\

5.  $V_2$ is convex with respect to $\varphi = Jx$ when $\xi$ is
fixed and with respect to $\xi$ when $x$ is fixed.\\

6.  $(||x||-||\xi||)^2 \leq V_2(Jx, \xi) \leq (||x||+||\xi||)^2 .$\\

7.  $V_2(Jx, \xi) \ge 0, \qquad \forall x,\xi \in B .$\\

8.  $V_2(Jx, \xi) = 0,$ only if  $ x=\xi .$\\

9.  $V_2(Jx, \xi) \rightarrow \infty,$  if  $ ||x||\rightarrow\infty$
(or
$||\xi||\rightarrow\infty$) and vice versa.\\

10. $V_2(Jx, \xi) \rightarrow 0 $ if $||x - \xi||\rightarrow 0, $ and
$||x||$,  $||\xi||$ are bounded, and vice versa.\\
\end {lem}

{\bf Proof} The continuty of the functional $V_2(Jx ,\xi)$ follows from
the continuty of norms and duality mappings in smooth spaces.

A differentiability of norms in uniformly smooth spaces gives the
differentiability  of $V_2(Jx ,\xi)$
with respect to $\varphi = Jx$ when $\xi$ is
fixed and with respect to $\xi$ when $x$ is fixed.

The formului for $ grad_{\varphi} V_2(Jx, \xi)$ and $ grad_{\xi} V_2(Jx, \xi)$
are verified by the direct calculations.

A monotonisity of operators $ grad_{\varphi} V_2(Jx, \xi)$ and
$ grad_{\xi} V_2(Jx, \xi)$ gives a convexity of the functional $V_2(Jx ,\xi)$.

The property 6 is obtained from expressions
$$||Jx||^2_{B^{*}} -2 <Jx, \xi> +||\xi||^2 \leq  ||x||^2 + 2||x||||\xi|| +
||\xi||^2 $$
and
$$||Jx||^2_{B^{*}} -2 <Jx, \xi> +||\xi||^2 \ge  ||x||^2 - 2||x||||\xi|| +
||\xi||^2 $$
It is obviously that the properties  7 - 9  follow from the property 6, and
the property 10 can be obtained from next statement as a corollary . Lemma
is proved.

There is a link between the functional $V_2(Jx ,\xi)$ and the traditional
functional $V_1(x ,\xi) = ||x - y||^2$ which is given by means of the
following statement.
\begin {lem} \label{s2}
The inequalities
$$2L^{-1} \delta_B (||x-\xi|| / 2C) \leq V_2(Jx, \xi) \leq
L^{-1} \rho_B (8LC||x-\xi||)$$
hold for all $x$ and $\xi$ in $B$,
where constant $L$ is from Theorem \ref{f54} and
$$C=2 \max \lbrace  1, \sqrt{(||x||^2 +||\xi||^2)/2} \rbrace.$$
\end {lem}
{\bf Proof} We are starting the left inequality. In
\cite{lz} it was shown that if $x \in B, y \in B $ and $||x||^2 + ||y||^2 = 2$,
then
$$||(x + y)/2||^2 \leq 1 - \delta_B (||x-y|| / 2)$$
where $\delta_B (\epsilon)$
is modulus of convexity of space $B$. It is known that in uniformly
convex Banach space, $\delta_B(\epsilon)$ is convex strictly
increasing function,
$0 \leq \delta_B (\epsilon) < 1,  \delta_B (0) = 0$.
We denote $R^2_1 = 2^{-1}
(||x||^2 + ||y||^2)$ ($x$ and $y$ are not zero at the same time) and introduce
new variables $ \tilde x = x/R_1$ and $ \tilde y = y/R_1$. Then
$$||\tilde x||^2 + ||\tilde y||^2 = R^{-2}_1(||x||^2 + ||y||^2) = 2 $$
Therefor, for $\tilde x$ and $\tilde y$ the inequality
$$||(\tilde x + \tilde y)/2||^2 \leq 1 -
\delta_B (||\tilde x - \tilde y||/2)$$
is valid.

Let us return to the old variables $x$ and $y.$ We obtain
$$||(x + y)/2R_1||^2 \leq 1 - \delta_B (||x-y|| / 2R_1)$$
or
\begin{equation} \label{a3}
||(x + y)/2||^2 \leq (||x||^2 + ||y||^2)/2 - R^2_1 \delta_B (||x-y|| / 2R_1),
\end{equation}
$$R_1 = \sqrt{(||x||^2 +||y||^2)/2} $$
Consider two cases:

1. $R_1 \ge 1.$ Then
\begin{equation} \label{a4}
||(x + y)/2||^2 \leq (||x||^2 + ||y||^2)/2 -  \delta_B (||x-y|| / 2R_1)
\end{equation}

2. $R_1 \leq 1.$ From (\ref{a3}) and the inequality  (see \cite{fg})
$$\epsilon ^2 \delta_B (\eta) \ge (4L)^{-1} \eta ^2 \delta_B (\epsilon),
\;\;\;\forall \eta \ge \epsilon >0 $$
we have
$$\delta_B (||x-y|| / 2R_1) \ge R^{-2}_1 (4L)^{-1} \delta_B (||x-y|| / 2) $$
and
\begin{equation} \label{a5}
||(x + y)/2||^2 \leq (||x||^2 + ||y||^2)/2 - (4L)^{-1} \delta_B (||x-y|| / 2)
\end{equation}
In vertue of $\min \lbrace 1, (4L)^{-1} \rbrace = (4L)^{-1},$
(\ref{a4}) and (\ref{a5}) if joined together give
$$||(x + y)/2||^2 \leq (||x||^2 + ||y||^2)/2 - (4L)^{-1}
\delta_B (||x-y|| / C_2),\;\; C_2 = 2\max \lbrace 1,R_1 \rbrace $$
If $||x|| \leq R$ and $||y|| \leq R$, then $C_2 = 2\max \lbrace 2,R \rbrace$
is the absolute constant.

Let $\varphi (x) = ||x||^2/2.$ Then
\begin{equation} \label{a6}
\varphi ((x +y)/2) \leq 2^{-1} \varphi (x)  + 2^{-1} \varphi (y) -
(8L)^{-1} \delta_B (||x-y|| / C_2),
\end{equation}
The formula  (\ref {a6}) shows (see Lemma \ref{f77}) that the functional
$\varphi (x)$ is uniformly convex
with modulus of convexity
$$\delta (||x-y||) = (4L)^{-1} \delta_B (||x-y|| / C_2). $$
Let us apply  Lemma \ref{f77} again. Because
$\varphi (x)$ is differentiable functional and
$\varphi'(x) = grad (||x||^2/2) = Jx$, then
$$<Jx - Jy, x - y> \ge (2L)^{-1} \delta_B (||x-y||/C_2)$$

Let $y$ be an arbitrary fixed point in $B$. In vertue of the uniform
convexity of $||x||^2/2$ we have
$$||(x + y)/2||^2 - ||x||^2 \ge <Jx,y - x> +
(2L)^{-1} \delta_B (||x-y||/2C_2).$$
Using now the identity $<Jx,x> = ||x||^2$, we obtain the following expression
$$<Jx,y> \leq ||(x + y)/2||^2 - L^{-1} \delta_B (||x-y||/2C_2).$$
So far as $g_B(\epsilon) / \epsilon$ is a nondecreasing
function, it is
$2 \delta_B (\epsilon/2) \leq \delta_B (\epsilon).$ Hence, one can write
$$ V_2(Jx, y) = ||x||^2 - 2<Jx,y> + ||y||^2 \ge$$
$$ ||x||^2 + ||y||^2 - 2^{-1}||(x + y)||^2 +
L^{-1} \delta_B (||x-y||/2C_2) \ge$$
$$(2L)^{-1}\delta_B (||x-y||/C_2) + L^{-1} \delta_B (||x-y||/2C_2) \ge$$
$$ 2L^{-1} \delta_B (||x-y||/2C_2).$$
The left estimate from Lemma is proved. Let us prove that right estimate
is valid.

Analogously to previous inequality, we can obtain the estimate
\begin{equation} \label{a7}
<Jx - Jy, x - y> \ge (2L)^{-1} \delta_{B^{*}} (||Jx-Jy||_{B^{*}}/C_2)
\end{equation}
for uniformly smooth space $B.$  From (\ref{a7}) one has
$$||Jx-Jy||_{B^{*}}||x - y|| \ge (2L)^{-1} \delta_{B^{*}}
(||Jx-Jy||_{B^{*}}/C_2)$$
Since $g_{B^{*}} (\epsilon)  = \delta_{B^{*}} (\epsilon)  / \epsilon $,
we can write
\begin{equation} \label{a8}
g_{B^{*}} (||Jx-Jy||_{B^{*}}/C_2) \leq C_4L||x - y||
\end{equation}
{}From geometry of Banach spaces  \cite{lz} it follows that
$$\rho _B (\tau) \ge \epsilon \rho /2 - \delta _{B^{*}} (\epsilon),\;\;
0 \leq \epsilon \leq 2,\;\; \tau > 0 $$
Therefore
$$\rho _B (4\delta _{B^{*}}(\epsilon)/ \epsilon) \ge
\delta _{B^{*}}(\epsilon)$$
We denote $h _B (\tau) = \rho _B (\tau) /\tau $. Then
$$h _B (4g_{B^{*}} (\epsilon)) \ge \epsilon / 4 $$
Setting
$$\epsilon = ||Jx-Jy||_{B^{*}}/C_2, $$
and using nonincreasing function $h _B (\tau)$, we obtain from (\ref{a8})
$$h _B (4g_{B^{*}} (\epsilon)) \leq h_B (8C_2L||x - y||) $$
Therefore
$$ ||Jx-Jy||_{B^{*}} \leq 4C_2 h_B (8C_2L||x - y||) $$
Using the inequality of Cauchy-Bunyakovski-Schwarz, we prove  the estimate
$$<Jx - Jy, x - y>  \leq (2L)^{-1} \rho _B (8C_2L||x - y||) $$
On the other hand, by the convexity of $V_2(Jx, y)$ one can obtain
$$V_2(Jx, y) \leq V_2(Jy, y) + 2<Jx - Jy, x - y> = 2<Jx - Jy, x - y>.$$
Finnaly, that implies the right estimate of Lemma.

\begin{re}
In the proving of Lemma \ref{s2},  the estimates for
normalized duality mappings in $B$
$$(2L)^{-1} \delta_B (||x-y||/C_2) \leq <Jx - Jy, x - y>
\leq (2L)^{-1} \rho _B (8C_2L||x - y||)$$
and
$$ ||Jx-Jy||_{B^{*}} \leq 4C_2 h_B (8C_2L||x - y||),\;\; 1<L<3.18,\;\;
h _B (\tau) = \rho _B (\tau) /\tau $$
\end{re}
were obtained (cf. \cite{aln1,aln9}). They are a quantitative description
of well known statement: in uniformly convex (uniformly smooth) Banach spaces
normalized duality the mapping $J$ is uniformly monotone
(uniformly continuous) operator on each bounded set.
\begin{re}
There is also a connection between the functional $V_2(Jx, \xi)$ and the
Young-Fenchel transformation, because
 $$||Jx||_{B^{*}} = \sup_{\xi\in B} \lbrace 2<J x,\xi>  -
||\xi||^2\rbrace .$$
\end{re}
\begin{re}
Notice  also that (\ref {k64}) is equivalent to
 $$V_2(Jx ,\xi) = ||x||^2 -2 <Jx, \xi> +||\xi||^2.$$
However, previous form (\ref {k64}) is necessarily used to obtain
properties 4 of the functional $V_2(Jx, \xi)$ and Lemma \ref{s2}.
\end{re}

Now we are ready to introduce first generalized projection operator in Banach
space.
\begin {de}  \label{gpo1}
Operator $\Pi_\Omega: B \rightarrow \Omega \subset B$ is called
the generalized projection operator if
it puts the arbitrary fixed point $ x \in B$ into the correspondence
with
the point of minimum for the functional $V_2(Jx, \xi)$ according
to the minimization problem
 $$\Pi_{\Omega}x = \hat x; \;\;\; \hat  x:  V_2(Jx ,\hat x) =
\inf _{\xi \in \Omega} V_2(Jx ,\xi) .$$
\end{de}
\begin{re}
In Hilbert space $ V_2(Jx ,\xi) = V_1(x ,\xi) $ and $\hat x =\bar x$.
\end{re}

 Next we describe the properties of the operator $ \Pi_\Omega$:\\

{\bf  7.a.} The operator $\Pi_{\Omega}$ is fixed in each point $\xi \in
\Omega $, i.e.
 $\Pi_{\Omega}\xi = \xi.$\\
{\bf Proof} follows from the property 8 of the functional $V_2(Jx, \xi)$)
because
$$ \inf _{\eta \in \Omega} V_2(J\xi ,\eta)  = V_2(J\xi ,\xi) = 0.$$

{\bf  7.b.}  $\Pi_{\Omega}$ is monotone ({\it d-accretive}) in $B$ , i.e.
 $$< Jx - Jy, \hat x -\hat y > \ge 0  .$$
{\bf Proof} By the Definition   \ref{gpo1}  of the operator
$\Pi_{\Omega},$ two inequalities
$$V_2(Jx_1 ,\xi) \ge   V_2(Jx_1 ,\hat x_1),$$
$$V_2(Jx_2 ,\eta) \ge   V_2(Jx_2 ,\hat x_2)$$
are satisfied for all $\xi,\eta \in \Omega $ and for all $ x,y \in B $.
Let be $\xi = \hat x_2$ and
$ \eta = \hat x_1$. Then
$$V_2(Jx_1 ,\hat x_2) \ge   V_2(Jx_1 ,\hat x_1),$$
$$V_2(Jx_2 ,\hat x_1) \ge   V_2(Jx_2 ,\hat x_2).$$
Then we have
$$ {||x_1||^2 -2 <Jx_1, \hat x_2> +||\hat x_2||^2}  +
{||x_2||^2 -2 <Jx_2, \hat x_1> +||\hat x_1||^2}  \ge  $$

$$ {||x_1||^2 -2 <Jx_1, \hat x_1> +||\hat x_1||^2}  +
{ ||x_2||^2 -2 <Jx_2, \hat x_2> +||\hat x_2||^2} .$$
{}From this one follows the expression
$$<Jx_1, \hat x_1> + <Jx_2, \hat x_2>{\quad } \ge{\quad } <Jx_1, \hat x_2> +
<Jx_2, \hat x_1>$$
and
$$  <Jx_1 - Jx_2, \hat x_1 - \hat x_2> \ge 0,{\qquad }
\forall x_1, x_2 \in B .$$

{\bf  7.c.}  $ < Jx - J\hat x, \hat x -\xi > \ge 0,  {\qquad }
\forall \xi \in\Omega  .$  \\
The point $\hat x$ is generalized projection of $x$  on $\Omega \subset B$
iff the inequality 7.c is satisfied.
{\bf Proof}   Definition   \ref{gpo1} gives
$$V_2(Jx ,\hat x) \leq   V_2(Jx ,\hat x + \theta (\xi - \hat x))$$
where $\theta \in [0,1] $ .
Using the property 5 of the convex functional $V_2(Jx ,\xi)$ we get
$$0 \ge V_2(Jx ,\hat x) -  V_2(Jx ,\hat x + \theta (\xi - \hat x)) \ge
 2 <J(\hat x + \theta (\xi - \hat x)) - Jx, \hat x - \hat x -
\theta (\xi - \hat x)> .$$
This leads to inequality
$$  <J(\hat x + \theta (\xi - \hat x)) - J x, \xi - \hat x> \ge 0. $$
Letting   $\theta \rightarrow 0,$  we have
$$  <J \hat x - J x,  \xi - \hat x> \ge 0,
{\qquad }  \forall \xi \in\Omega. $$
Inversely, if 7.c is filfilled then in vertue of the property 4 of
$V_2(Jx ,\xi) $ one occurs
$$ V_2(Jx ,\xi) - V_2(Jx ,\hat x) \ge 2<J \hat x - J x,  \xi - \hat x> \ge 0,
{\qquad }  \forall \xi \in\Omega, $$
i.e.
$$ V_2(Jx ,\xi) \ge V_2(Jx ,\hat x), {\qquad }  \forall \xi \in\Omega. $$
Therefore  $\hat x = \Pi_{\Omega}x .$
\begin{re}
The property 7.c we call {\it basic variational principle}
for $\Pi_{\Omega}$ in $B$ since it is the eqivalent representation
of the solution (in a variational inequality form) for minimization problem
of the Lyapunov functional  $V_2(Jx, \xi),$ (cf. 4.c, 5.c).
\end{re}
{\bf  7.d.}  $ < Jx - J\xi, \hat x -\xi > \ge 0,  {\qquad }
\forall \xi \in\Omega .$ \\
{\bf Proof} It is obviously that
$$  V_2(Jx,\hat x) \leq V_2(Jx,\xi) +  V_2(J \xi,\hat x). $$
Rewriting this inequality in detail acording to equalities
$$ V_2(Jx,\hat x) = ||x||^2 -2 <Jx, \hat x> +||\hat x||^2, $$
$$V_2(Jx,\xi) =  ||x||^2 -2 <Jx, \xi> +||\xi||^2, $$
$$ V_2(J \xi,\hat x) = ||\xi||^2 -2 <J\xi, \hat x> + ||\hat x||^2 $$
we obtain
$$||x||^2 -2 <Jx, \hat x> + ||\hat x||^2 \leq ||x||^2 -2 <Jx, \xi> +||\xi||^2 +
  ||\xi||^2 -2 <J\xi, \hat x> +||\hat x||^2. $$
{}From this it follows that
$$<Jx, \xi> +  <J\xi, \hat x> -  <Jx, \hat x> \leq ||\xi||^2 .$$
We have now the property 7.d. because $ ||\xi||^2 = <J\xi, \xi> $.

\begin{re}
{}From inequality $V_2(Jx,\hat x) \leq V_2(Jx,\xi)$, the following useful
estimate
$$ 2<Jx, \xi> \leq ||\hat x||^2 + ||\xi||^2 $$
can be obtained.
\end{re}

{\bf  7.e. }  $ < Jx - J\hat x ,  x -\xi > \ge 0,  {\qquad }
\forall \xi \in\Omega .$ \\
{\bf Proof}  Using basic variational principle for $\Pi_{\Omega}$ in $B$
(the property 7.c) we can write
$$ < Jx - J\hat x, \hat x - x > +  < Jx - J\hat x, x -\xi >  \ge 0 .$$
Acording to the property 7.b
$$ < Jx - J\hat x, \hat x - x > \leq 0 .$$
Therefore the inequality 7.e is valid.

{\bf  7.f.} $||\hat x - \hat y|| \leq C g^{-1}_B (2LC^2 g^{-1}_{B^{*}}
(2LC ||x-y||))$,
where constant $L$ is from the Theorem \ref{f54} and
$$C = 2 \max \lbrace  1, ||x||,||y||,||\hat x||,||\hat y|| \rbrace .$$

{\bf Proof}  From  \cite{aln9,aln1} it follows that
$$(2L)^{-1}\delta_B (||\hat x_1 - \hat x_2||/C_1) \leq
<J \hat x_1 - J \hat x_2, \hat x_1 - \hat x_2>$$
where
$$C_1(||\hat x_1||, ||\hat x_2||) =
 2 \max \lbrace  1, \sqrt{(||\hat x_1||^2
+||\hat x_2||^2)/2}\rbrace .$$
Using the propertie 7.c of the projection operator $\Pi_\Omega$ we obtain
$$ <J \hat x_1 - J \hat x_2, \hat x_1 - \hat x_2> =
   <J x_1 - J \hat x_2, \hat x_1 - \hat x_2> -
   <J x_1 - J \hat x_1, \hat x_1 - \hat x_2> \leq $$
$$ <J x_1 - J \hat x_2, \hat x_1 - \hat x_2> .$$
Using again the propertie 7.c one can write
$$ <J x_1 - J \hat x_2, \hat x_1 - \hat x_2> \leq
   <J x_1 - J \hat x_2, \hat x_1 - \hat x_2> -
   <J x_2 - J \hat x_2, \hat x_1 - \hat x_2> = $$
$$ <J x_1 - J x_2, \hat x_1 - \hat x_2> .$$
It is obviously that
$$ <J x_1 - J x_2, \hat x_1 - \hat x_2> \leq ||J x_1 - J x_2||_{B^{*}}
||\hat x_1 - \hat x_2|| .$$
Therefore  the inequalities
$$(2L)^{-1}\delta_B (||\hat x_1 - \hat x_2||/C_1) \leq
||J x_1 - J x_2||_{B^{*}} ||\hat x_1 - \hat x_2||$$
and
$$||\hat x_1 - \hat x_2|| \leq C_1 g^{-1}_{B} (2C_1 L||J x_1 - J
x_2||_{B^{*}})$
are valid. From \cite{aln9,aln1,aln2} we have
$$||J x_1 - J x_2||_{B^{*}} \leq C_2 g^{-1}_{B^{*}}
{(2C_2 L||\hat x_1 - \hat x_2||)}$$
where
$$C_2(||x_1||, ||x_2||) =
 2 \max \lbrace  1, \sqrt{(||x_1||^2
+||x_2||^2)/2}\rbrace $$
This leads to final result 7.f.
\begin{re}   If  $||x ||\leq R, || \hat x || \leq R,  ||y ||\leq R $ and
 $|| \hat y|| \leq R,$ then  $C = 2 {max} \lbrace 1, R \rbrace$  is
absolute constant and 7.f expresses the uniform continuity  of operator
 $\Pi_{\Omega} $ in Banach space on each bounded set.
\end{re}

{\bf  7.g.}   $ < Jx - Jy, \hat x -\hat y > \ge (2L )^{-1}\delta_B
(||\hat x -\hat y||/C) ,$  where
$$C = 2 \max \lbrace  1, ||\hat x||, ||\hat y|| \rbrace .$$
This inequality contains in the proof of the property 7.f.

{\bf  7.h.}  The operator $ \Pi_{\Omega}$ gives absolutely best
approximation of
$x\in B$ with respect to
 functional $V_2(Jx, \xi)$, i.e.
 $$V_2(J \hat x, \xi) \leq V_2(Jx, \xi) - V_2(Jx, \hat x), {\qquad }
 \forall \xi \in\Omega .$$
Consequently, $ \Pi_{\Omega}$ is {\it conditionally nonexpanseve \/}
with respect to functional $V_2(Jx, \xi)$
operator in Banach space, i.e.
 $$V_2(J \hat x, \xi) \leq V_2(Jx, \xi), {\qquad } \forall \xi \in\Omega .$$
{\bf Proof} Using the property 7.c we obtain
$$ < Jx , \hat x - \xi > \leq < J \hat x, \hat x - \xi >,  {\qquad }
\forall \xi \in\Omega .$$
{}From this it follows that
$$ < Jx , \xi > \ge < J \hat x, \xi > + < Jx , \hat x > -
< J \hat x, \hat x> = $$
$$ < J \hat x, \xi > + < Jx , \hat x > - ||\hat x||^2 $$
It is equivalent to the inequality
$$||x||^2 -2 <Jx, \xi> + ||\xi||^2 \leq $$
$$||\hat x||^2 -2 <J\hat x, \xi> + ||\xi||^2 + ||x||^2 - 2 <Jx, \hat x> +
||\hat x||^2. $$
The property 7.h is derived now in vertue of following expressions
$$ V_2(Jx,\hat x) = ||x||^2 - 2 <Jx, \hat x> +||\hat x||^2, $$
$$V_2(Jx,\xi) =  ||x||^2 -2 <Jx, \xi> +||\xi||^2, $$
$$ V_2(J \hat x,\xi) = ||\hat x||^2 -2 <J \hat x, \xi> + ||\xi||^2 .$$

{\bf 7.i.} Any $\Pi_{\Omega}$ satisfies the inequality
$$<(J - J\Pi_{\Omega})x - (J - J\Pi_{\Omega})y, \Pi_{\Omega}x -
\Pi_{\Omega}y> \ge 0,
\;\;\forall x,y \in B.$$

{\bf Proof} From the property 7.c it follows that
$$< Jx - J\hat x, \hat x - \hat y > \ge 0$$
and
$$< Jy - J\hat y, \hat y - \hat x > \ge 0$$
It gives in sum 7.i.

\begin{re}  If B=H, then the formulas 7.a - 7.e and 7.h - 7.i coincide with
ones
4.a - 4.e and 4.h - 4.i, but 7.f  and 7.g differ from 4.f and 4.g by only
constants
(on any bounded set).
 \end{re}
 Using properties of  the generalized projection operator $\Pi_\Omega$
we
obtained the  theorem.
\begin {th} \label{g4}
The following assertions hold for the alternating method in Banach space
(method of successive generalized projections):
(\ref{f16}) and (\ref{f17}):\\
\\
1)  $V_2(Jx _{n+1},\xi) \leq  V_2(Jx _{n},\xi), {\qquad } \forall \xi
\in\Omega_{*} = \bigcap _{i=1}^{m} \Omega_{i} .$\\
\\
2) There exists  a subsequence $\lbrace x_{n_k} \rbrace $ of the
sequance $ \lbrace x_{n}
 \rbrace$ such that  $x_{n_k} \rightarrow x_{*} $  weakly,
 where $ x_{*}\in  \Omega _{*} .$\\
\\
3) If $J$ is sequential weakly continuous operator then
$x_{n} \rightarrow x_{*} $ weakly.\\
\\
If $\lbrace x_{n} \rbrace $ is an ordered sequence of the elements
$x^{j}_i,\;i = 0,1,...,\;j = m,m-1,...,1,$ such that $x^{m}_i =
\Pi_{m}{x^{1}_{i-1}};\; x^{j}_i = \Pi_{j}{x^{j+1}_i},\;j = m-1,m-2,...2,1;\;
x^{1}_{-1} = x_0 $, then, in addition to 1) - 3),\\
\\
4)  $\sum _{n=0}^{\infty} V_2(Jx _{n}, x_{n+1}) < \infty .$\\
\\
5)  $ ||x_{n} - x_{n+1}|| \rightarrow 0 , \;\;$ for $n \rightarrow \infty .$\\
\end {th}
{\bf Proof} Let $\xi \in \Omega $ be an arbitrary  fixed point and
let $x_n = x^{(1)}_n \in \Omega_1$ is $n$-approximation, determined by
(\ref{f16}) and
(\ref{f17}). We denote ${{\hat x} _n}^{(1),m} = \Pi _{\Omega _m} x^{(1)}_n $.
Then in vertue of 7.h  we obtain
$$V_2(J {{\hat x} _n}^{(1),m}, \xi) \leq V_2(Jx^{(1)} _n, \xi). $$
Now we denote  ${{\hat x} _n}^{(1),m,m-1} = \Pi _{\Omega _{m-1}}
{{\hat x} _n}^{(1),m} $.
Again, in vertue of 7.h we can write
$$V_2(J {{\hat x} _n}^{(1),m,m-1}, \xi)
\leq V_2(J {{\hat x} _n}^{(1),m}, \xi). $$
Continuiting this process further, we have at finish
$ {{\hat x} _n}^{(1),m,m-1,...,2,1} = \Pi _{\Omega _1} {{\hat x} _n}^
{(1),m,m-1,...,2}$ and
$$V_2(J {{\hat x} _n}^{(1),m,m-1,...,3,2,1}, \xi)
\leq V_2(J {{\hat x} _n}^{(1),m,m-1,...,3,2}, \xi). $$
These inequalities give the estimate
$$V_2(J {{\hat x} _n}^{(1),m,m-1,...,2,1}, \xi) \leq V_2(Jx^{(1)}_n, \xi). $$
Because
$${{\hat x} _n}^{(1),m,m-1,...,2,1} =
\Pi _{\Omega _1} \circ \Pi _{\Omega _2}  \;\circ \cdotp \cdotp \cdotp \circ \;
\Pi _{\Omega _m} x^{(1)}_n  = x^{(1)}_{n+1} = x_{n+1},$$
the assertion 1) is hold for $n = 1,2,...$.\\
Now we conclude by the same way that
$$V_2(Jx^{(1)}_{1},\xi) \leq V_2(Jx_{0},\xi), \;\;x^{(1)}_{1} = \Pi x_{0}.$$
Thus, nonnegative sequence $V_2(Jx^{(1)}_{n},\xi)$ is bounded from above by
$V_2(Jx_{0},\xi).$ Therefore it
has a limit. Further, since the  property 9 for the functional $V_2(Jx,\xi)$
holds, we obtain the estimate
$$  ||x^{(1)}_{n}|| \leq ||x _{0}|| + 2 ||\xi|| = const,\;\; n=1,2,....$$
In reflexive Banach space any bounded points set is weakly compact. Then
it follows that one can select the weakly convergent to
$x^{*} $ subsequence
$\lbrace x^{(1)}_{n_k} \rbrace $ from the sequence
$\lbrace x ^{(1)}_{n} \rbrace $.

Let us show that $x^{*} \in \Omega_{*}.$   In fact, $x^{*} \in \Omega_{1}$
becau
all of the sets $\Omega_{i}$ are closed (i.e. weakly closed).
Consider now the siquence  $x^{(2)}_{n}$ such that
$$x^{(2)}_{n+1} = \Pi^{n+1}x^{(2)}_0,  {\qquad } n=0,1,2,...,$$
where $x ^{(2)}_0 =  \Pi_{2} \circ \Pi_{3} \;\circ \cdotp \cdotp \cdotp
\circ \;\Pi_{m} {x_0} $ and
$$\Pi =  \Pi_{2} \circ \Pi_{3} \;\circ \cdotp \cdotp \cdotp \circ \;
 \Pi_{m}\;\circ\;\Pi_{1},  \qquad      \Pi_{i} = \Pi_{\Omega_{i}}.$$
It is obviously that $x^{(2)}_{n} \in \Omega_{2}.$  We choose subsequence
$\lbrace x^{(2)}_{n_k} \rbrace \in \lbrace x^{(2)}_{n} \rbrace$ for which
$ \Pi_{\Omega_{1}} x^{(2)}_{n_k} = x^{(1)}_{n_k}.$
In this case 7.h gives
$$V_2(J x^{(2)}_{n_k}, x^{(1)}_{n_k}) \leq V_2(Jx^{(2)}_{n_k}, \xi)
- V_2(Jx^{(1)}_{n_k}, \xi).$$
We know that $V_2(Jx^{(2)}_{n_k}, \xi)$  and $V_2(Jx^{(1)}_{n_k}, \xi)$ have
limits therefore
$$ V_2(Jx^{(2)}_{n_k}, \xi) -  V_2(Jx^{(1)}_{n_k}, \xi) \rightarrow 0 ,
\;\; as\;\;  n_{k} \rightarrow \infty .$$
Using the property 10 of functional $V_2(Jx, \xi)$ we obtain
$$V_2(J x^{(2)}_{n_k}, x^{(1)}_{n_k}) \ge 2L^{-1} \delta_B
(||x^{(2)}_{n_k} - x^{(1)}_{n_k}|| / C) \rightarrow 0 ,$$
because $C$ is absolute constant. This means that
$$||x^{(2)}_{n_k} - x^{(1)} _{n_k}|| \rightarrow 0,\;\; as\;\; n_k \rightarrow
\infty ,$$
$x^{(2)} _{n_k}$  converges weakly to $x^{*}$ and $x^{*} \in \Omega_{2}.$
By the same way we can show that $x^{*} \in \Omega_{3},\; x^{*} \in \Omega_{4},
..., x^{*} \in \Omega_{m}.$
Thus, $ x_{*} \in  \Omega _{*}= \bigcap _{i=1}^{m} \Omega_{i},$ and we proved
the the assertion 2). \\
Let us show now that the whole siquence $\lbrace x _{n} \rbrace $
converges to $x^{*} $ weakly under the condition 3).
Suppose that there exists $x^{**}$ such that
$x_n \rightharpoonup x^{**}.$ We can write
 $$V_2(Jx_n ,x^{*}) = ||x_n||^2 -2 <Jx_n, x^{*}> + ||x^{*}||^2.$$
 $$V_2(Jx_n ,x^{**}) = ||x_n||^2 -2 <Jx_n, x^{**}> + ||x^{**}||^2.$$
The siquences  $V_2(Jx_n ,x^{*})$ and $V_2(Jx_n ,x^{**})$ have limits,
therefore their difference has common  limit . Denote
$$a = \lim_{n \rightarrow \infty} (V_2(Jx_n ,x^{*}) - V_2(Jx_n ,x^{**}))$$
It is easy to see that
$$a = 2 \lim_{n \rightarrow \infty} <Jx_n, x^{**} - x^{*}> + ||x^{*}||^2 -
||x^{**}||^2 .$$
If $x_n \rightharpoonup x^{*}$, then in vertue of
sequential weak continuty of operator $J$, we have
$$a = 2<Jx^{*}, x^{**} - x^{*}> + ||x^{*}||^2 - ||x^{**}||^2.$$
If $x_n \rightharpoonup x^{**}$, then
$$a = 2<Jx^{**}, x^{**} - x^{*}> + ||x^{*}||^2 - ||x^{**}||^2.$$
{}From this it follows that
$$ 0 =  <Jx^{**} - Jx^{*}, x^{**} - x^{*}> \ge L^{-1} \delta_B
(||x^{**} - x^{*}||/C),$$
where $C$ is absolute constant. Therefore $x^{**} = x^{*}.$

Let us prove now 4) and 5). For the ordered sequence of the elements
$\lbrace x_{n} \rbrace $ = $\lbrace x^{j}_i \rbrace $ we have
$$ x_{0} = x_{0}, x_{1} = x^{m}_{0}: \;\;V_2(J x^{m}_{0}, \xi)
\leq V_2(Jx_{0}, \xi) - V_2(Jx_{0}, x^{m}_{0}).$$
$$ x_{1} = x^{m}_{0}, x_{2} = x^{m-1}_{0} : \;\;V_2(J x^{m-1}_{0}, \xi)
\leq V_2(Jx^{m}_{0}, \xi) - V_2(Jx^{m}_{0}, x^{m-1}_{0}).$$
$$ x_{n}, x_{n+1}: \;\;V_2(J x_{n+1}, \xi)
\leq V_2(Jx_{n}, \xi) - V_2(Jx_{n}, x_{n+1}).$$
Then
$$V_2(J x_{n+1}, \xi) \leq V_2(Jx_{0}, \xi) -  \sum _{i=0}^{n}
V_2(Jx_{i}, x_{i+1}) $$
{}From this it follows 4) because $V_2(J x_{n}, \xi)$ is nonnegative and it has
limit as $n \rightarrow \infty $. Therefore
$$V_2(J x_{n}, \xi) \rightarrow 0, \;\;n \rightarrow \infty  $$
The statement 5) is obtained now from Lemma \ref{s2} (the propety 10  of the
functional  $V_2(Jx, \xi))$. Theorem is completely proved.

In what follows  we discuss statement 3 from the theorem.
We recall that  $F$ is called a sequential weakly continuous mapping if
from the relation  $ x_{n} \rightarrow x $ (weakly)  it follows  that
$ Fx_{n} \rightarrow Fx $  (also weakly).

Theorem \ref{g4} is valid for dual mapping $ J^{\mu}$ with the gauge
function
$\mu (t) = t^{p-1}$, defined by the relations
$$||J^{\mu} x||_{B^{*}} = ||x||^{p-1},  {\qquad } <J^{\mu} x, x> = ||x||^{p},$$
i.e.
$$ J^{\mu} x = grad ||x||^{p}/p . $$
(Notice that normalized dual mapping  corresponds to $p=2$). We set
$$V_3(J^{\mu} x,\xi) = q^{-1}||J^{\mu} x||^q_{B^{*}}- <J^{\mu} x,\xi>  +
p^{-1}||\xi||^p,$$
$$ p^{-1} + q^{-1} = 1 .$$
The function $ q^{-1}||J^{\mu} x||^q_{B^{*}} $ is conjugate to the function
$ p^{-1}||x||^p $, i.e.
$$q^{-1}||J^{\mu} x||^q_{B^{*}} = \sup_{\xi\in B} \lbrace <J^{\mu} x,\xi>
- p^{-1}||\xi||^p \rbrace .$$
Therefore  $ V_3(J^{\mu} x,\xi) \ge 0 ,{\quad} \forall  x,\xi\in B$ . If  we
now define
$$ \hat {\hat x}: {\qquad } V_3(J^{\mu} x, \hat {\hat x}) =
\inf_{\xi\in\Omega}V_3(J^{\mu} x,\xi) $$
then it can be shown that\\
$$<J^{\mu} x - J^{\mu} \hat{\hat x}, \hat{\hat x} - \xi> \ge 0, {\qquad}
\forall\xi\in\Omega $$
and
$$V_3(J^{\mu} x, \hat{\hat x})\leq V_3(J^{\mu} x,\xi) - V_3((J^{\mu} \hat{\hat
x},\xi). $$
The last inequalities are used mainly in the proof of Theorem \ref{g4}
which will
be given in forthcoming paper.
The property of the uniform continuity of the operator $J^{\mu} x$ can be
obtained using
the results from \cite{not,zr}.
\begin{co}
 In Banach space $l^p, p>1$ the sequence $ x_{n}$ weakly convergse
 to\\  $ x_{*}\in  \Omega _{*}= \bigcap _{i=1}^{m} \Omega_{i}.$\\
 \end{co}
This holds because in  $l^p, \infty >p>1$, dual mapping $ J^{\mu}$ with the
gauge
function $\mu (t) = t^p$  is sequential weakly continuous  \cite{br}.
\begin {re}   $V_3(J^{\mu} x, \xi)$ and $V_2(Jx,\xi)$ coinside for $ p=2$
 (up to constant 2).
 \end{re}
It can be shown in a way
similar to the case of  metric projection operator $P_\Omega$ in Banach
space
\cite{aln2} that generalized projection operator $\Pi_\Omega$ is stable
with
respect to peturbation of the set  $\Omega$.

Let  $\Omega_1$ and  $\Omega_2$ be convex closed sets, $x \in B$ and
 $H (\Omega_1, \Omega_2) \leq \sigma ,$ where
$$ H(\Omega_1, \Omega_2) =
\max \lbrace \sup_{z_1\in \Omega_1} \inf _{z_2\in \Omega_2}
 ||z_1 - z_2||, {\quad}\sup_{z_1\in \Omega_2} \inf _{z_2\in
\Omega_1}||z_1 - z_2||\rbrace $$
is a Hausdorff distance between $\Omega_1$ and $\Omega_2.$
Let also $\Pi_{\Omega_1}x = \hat x_1, \Pi_{\Omega_2}x = \hat x_2$

\begin {th} \label{g6}
If $B$ is uniformly convex  Banach space, $\delta _B (\epsilon) $  is
modulus of the
convexity of  $B$ and $\delta^{-1}_B (\cdot)$ is an inverse function,
then
$$||\hat x _1 -\hat x_2|| \leq C_1\delta^{-1}_B (4LC_2\sigma ) ,$$
$$ C_1 = 2 \max \lbrace 1,||\hat x_1||, ||\hat x_2|| \rbrace ,$$
$$ C_2 = 2 \max \lbrace  ||Jx - J\hat x_1||_{B^{*}},
||Jx - J\hat x_2||_{B^{*}} \rbrace.$$
\end {th}
{\bf Proof}  We have
$$ <J\hat x_2 - J\hat x_1, \hat x_2 - \hat x_1> \ge (2L )^{-1}\delta_B
(||\hat x_1 -\hat x_2||/C_1) $$
Using the condition $H (\Omega_1, \Omega_2) \leq \sigma ,$ we obtain the
following: there exists the element $\xi_1 \in \Omega_1$ such that
$||\hat x_2 - \xi_1|| \leq \sigma $ and for all $x \in B$
$$<Jx - J \hat x_1, \hat x_2 - \hat x_1> =
<Jx - J \hat x_1, \hat x_2 - \xi_1> + $$
$$+ <Jx - J \hat x_1, \xi_1 - \hat x_1>
\leq \sigma ||Jx - J\hat x_1 ||_{B^{*}},$$
because
$$< Jx - J\hat x_1, \xi_1 -\hat x_1> \leq 0$$
Similarly to previous inequalities, we obtain for the element $\xi_2
\in \Omega_2$
$$<Jx - J \hat x_2, \hat x_2 - \hat x_1> =
<Jx - J \hat x_2, \hat x_2 - \xi_2> + $$
$$+ <Jx - J \hat x_1, \xi_2 - \hat x_1>
\leq \sigma ||Jx - J\hat x_2 ||_{B^{*}},$$
because
$$< Jx - J\hat x_2, \xi_2 -\hat x_1> \leq 0.$$
Therefore,
$$ <J\hat x_2 - J\hat x_1, \hat x_2 - \hat x_1> \leq
\sigma (||Jx - J\hat x_1 ||_{B^{*}} + ||Jx - J\hat x_2 ||_{B^{*}})$$
$$\leq 2 \sigma  \max \lbrace ||Jx - J\hat x_1 ||_{B^{*}},
||Jx - J\hat x_2 ||_{B^{*}} \rbrace = \sigma C_2$$
hold. The final result follows from expression
$$\sigma C_2 \ge  (2L )^{-1}\delta_B (||\hat x_1 -\hat x_2||/C_1) .$$
Theorem is proved.

If  $||x|| \leq R , ||\hat x_1|| \leq R $ and $||\hat x_2 || \leq R ,$
then $ C_1$
and $C_2$ are absolute constants, because operator $J$ is bounded in any
Banach space.

 \section{Generalized Projection Operator $\pi_\Omega$ in Banach Space}
\setcounter{equation}{0}
Here we introduce  generalized projection operator $\pi_\Omega$
in Banach space and describe its properties. Then  we use this operator
to establish equivalence between  solution of the variational inequality
in Banach space and solution of the corresponding operator equation.
In other words, we solve first problem described in \S 2.  Finally,
we obtain a link between operators $\pi_\Omega$ and $\Pi_\Omega$
by means of the normalized duality mappings $J$ and $J^{*}$.

We assume that $ \varphi$ is an arbitrary element of the space $ B^{*}$.
\begin {de}
 The generalized projection $ \tilde \varphi$ of the element  $ \varphi$
on the set $\Omega \subset B$ is given by means of a minimization
problem
 $$\pi_{\Omega}\varphi = \tilde \varphi ;\;\;\;    \tilde \varphi:
V_4(\varphi , \tilde \varphi) = \inf _{\xi \in \Omega} V_4(\varphi ,\xi)
$$
where
\begin{equation} \label{f72}
V_4(\varphi ,\xi) = ||\varphi||^2_{B^{*}} -2 <\varphi, \xi>
+||\xi||^2_{B}.
\end{equation}
 \end{de}
\begin {re}   In Hilbert space
$ V_4(\varphi ,\xi) = V_3((J^p \hat{\hat x},\xi) =
 V_2(Jx ,\xi) = V_1(x ,\xi) $ and\\
 $\varphi = \hat{\hat x }= \hat x =\bar x  $ .
 \end{re}
In what follows we list properties of the generalized projection
operator
$\pi_\Omega$ in Banach Space. Sinse the proofs are similar to Section 7,
we omit they.\\[2mm]

{\bf  8.a.}  The operator $\pi_{\Omega}$ is J-fixed in each point $\xi
\in \Omega$, i.e.
  $\pi_{\Omega}J\xi = \xi .$

{\bf  8.b.}  $\pi_{\Omega}$ is monotone in $B^{*}$, i.e. for all
$\varphi_1,\varphi_2 \in B^{*}$
 $$<\varphi_1 - \varphi_2, \tilde \varphi_1 - \tilde \varphi_2 > \ge 0 .
$$

{\bf  8.c.}  $ <\varphi - J \tilde\varphi, \tilde \varphi - \xi > \ge 0,
 {\qquad }        \forall \xi \in\Omega .$  \\
The point $\tilde\varphi$ is generalized projection of $\varphi$  on
$\Omega \subset B$
iff the inequality 8.c is satisfied.
The property 8.c we call {\it basic variational principle} for $\pi_{\Omega}$
in dual couple $B,B^{*}$.

{\bf  8.d.} $ <\varphi - J \xi, \tilde \varphi - \xi > \ge 0,  {\qquad }
\forall \xi \in\Omega .$

{\bf  8.e.} $ <Jx - J\tilde x, x - \xi> \ge 0$, where $ \tilde x =
\pi_{\Omega}Jx,
           {\qquad }        \forall \xi \in\Omega .$\\

{\bf  8.f.}  $||\tilde \varphi_1 - \tilde \varphi_2||
           \leq C g^{-1}_B (2LC ||\varphi_1 - \varphi_2||_B{*}), $
where the constant $L$ is from Theorem \ref{f54} and
$$C = 2 \max \lbrace  1, ||\tilde \varphi _1||,
||\tilde \varphi _2|| \rbrace .$$

{\bf  8.g.}  $ <\varphi_1 - \varphi_2, \tilde \varphi_1 - \tilde
\varphi_2 >
           \ge (2L )^{-1}\delta_B (||\tilde \varphi_1 - \tilde
\varphi_2||/C) ,$
where $C$ is constant from 8.f.
\begin {re} If  $ ||\varphi _1|| \leq R$, $||\varphi _2|| \leq R$,
$ ||\tilde \varphi _1|| \leq R   $, $||\tilde \varphi _2|| \leq R$ then
 $C = 2 \max \lbrace 1, R \rbrace $ is an absolute constant and  8.f
expresses
uniform continuity  of  the operator $ \pi_{\Omega} $ in Banach space on
each
 bounded set.
  \end{re}

{\bf  8.h.} The operator $ \pi_{\Omega}$ gives absolutely best
approximation of $x\in B$
with respect to functional $V_4(\varphi ,\xi) $
$$V_4(J\tilde \varphi,\xi) \leq  V_4(\varphi ,\xi) -  V_4(\varphi ,
\tilde \varphi). $$
Consequently, $ \pi_{\Omega}$ is {\it conditionally nonexpanseve \/}
with respect to functional $V_4(\varphi ,\xi) $
operator in Banach space, i.e.
$$V_4(J\tilde \varphi,\xi) \leq  V_4(\varphi ,\xi) .$$

{\bf 8.i.} Any $\pi_{\Omega}$ satisfies the inequality
$$<(I_{B^{*}} - J\pi_{\Omega})\varphi_1 -
(I_{B^{*}} - J\pi_{\Omega})\varphi_2,
\pi_{\Omega} \varphi_1 - \pi_{\Omega} \varphi_2 > \ge 0,
\;\;\forall \varphi_1, \varphi_2 \in B^{*},$$
where $I_{B^{*}} : B^{*} \rightarrow B^{*}$ is identical operator in $B^{*}$.

\begin {re}
Similarly to 7.i, the proof follows from the inequalities
$$ <\varphi_1 - J \tilde\varphi_1, \tilde \varphi_1 - \tilde\varphi_2> \ge 0 $$
and
$$ <\varphi_2 - J \tilde\varphi_2, \tilde \varphi_2 - \tilde\varphi_1> \ge 0 $$

Similarly to operator $\Pi_\Omega $, the generalized projection operator
$\pi_\Omega$  in Banach Space is stable with respect to
peturbation of the set  $\Omega$.
\end{re}

Using  the properties of generalized projection operator $\pi_\Omega$
we obtained
 Banach space analogue of Theorem \ref{th1}.

\begin {th} \label{f80}
 Let $ A $ be an arbitrary  operator from Banach space $B $ to $B^{*}$,
$\alpha$ an
 arbitrary fixed positive number, $f\in B^{*}$.  Then  the point $x \in
\Omega  \subset B$
is a solution of variational inequality
\begin{equation} \label {f982}
<Ax-f, \xi - x> \ge 0, {\qquad }        \forall \xi \in\Omega ,
\end{equation}
if and only if x is a solution of the operator equation in B
\begin{equation} \label {f482}
 x=\pi_{\Omega}(Jx - \alpha (Ax-f)).
\end{equation}
\end{th}
The proof follows from  the property 8.c  and eqivalent representation
$$<Jx - \alpha (Ax-f) - Jx, x - \xi> \ge 0,
{\qquad }  \forall \xi \in\Omega (x)$$
for variational inequality (\ref {f982}).

It is not hard to verify that
 $$\Pi _\Omega = \pi _\Omega J ,   {\qquad }   \pi _\Omega = \Pi _\Omega
J^{*} $$
 where $J: B\rightarrow B^{\ast}$ is a normalized duality mapping in $B$
 and $ J^{*}: B^{\ast} \rightarrow B $ is normalized duality mapping
in $  B^{\ast}$.  Therefore  we can rewrite (\ref{f482}) in the form of
$$x=\Pi_{\Omega}J^{*}(Jx - \alpha (Ax-f)) .$$
Denote also that $ J^{*} = J^{-1}. $

\begin {re}
For the iterative process (\ref{k6}) we proved the assertions being
analogous to the
Theorems 3 and 4 from the paper \cite{aln1} (see also Remark  7  in
\cite{aln1}).
\end {re}

It is interesting to note that $\pi_\Omega = J^{*}$ in the case $\Omega
= B$
(the problem of solving the equation $Ax = f).$ Then (\ref{f482}) is
rewritten as
$$x=J^{*}(Jx - \alpha (Ax-f)) $$
or in the form of
\begin{equation} \label {f84}
 Jx = Jx - \alpha (Ax-f)
\end{equation}
because $JJ^{*} = I_{B^{*}}$. Here $I_{B^{*}}:B^{\ast} \rightarrow B^{\ast}$
is identical operator. Iterative method for (\ref{f84})
$$Jx_{n+1} = Jx_{n} - \alpha_n (Ax _{n}-f),{\quad }     n=0,1,2,...,
x_o \in B, {\quad }  x_n=J^{*}Jx_n $$
has been studied earlier in  \cite{aln1,aln8}.

Along with (\ref{k6}) we considered the following iterative processes:
$$x_{n + 1} = \pi_{\Omega}(Jx_{n}-\alpha_n (Ax_{n} - f) /||Ax_{n} -
f||),
{\qquad }n=0,1,2... $$
for the variational inequality (\ref{f5}) with nonsmooth unbounded
operator $A,$
and
$$x_{n + 1} = \pi_{\Omega}(Jx_{n}-\alpha_n ( \partial
u(x_{n})/||\partial u(x_{n})||),
{\qquad }n=0,1,2... $$
and
$$x_{n + 1} = \pi_{\Omega}(Jx_{n}-\alpha_n ( u(x_{n}) - u^{*})\partial
u(x_{n})/
||\partial u(x_{n})||^2),
{\qquad }n=0,1,2... $$
for the minimization of functional $u(x).$ Here $u^{*} =
\min_{x  \in  \Omega} u(x).$

Finally, more general variational inequalities in Banach spaces can be
considered. (See (\cite{jy}).
\begin {re}
In  Section 5 and Section 8 we formulated new equivalence theorems between
variational
inequalities in Banach spaces and corresponding operator equations
(\ref{f282}) and (\ref{f482}) with
metric projection operator $P_{\Omega}$  and  generalized projection
operator $\pi_\Omega$, respectively. It is natural to call the equations
of this type $direct$ $projection$ $equations$. At the same time one can be
established the connection between
variational inequalities and other operator equations, so called
Wiener-Hopf equations, including  quasivariational
inequalities and complementarity problems. (See, for example,
\cite{jy,no,shi}).
\end {re}
In conclusion , we notice  that generalized projection
operators in Banach spaces constructed above and metric projection
operators in
Hilbert space are defined in similar way in the form of minimization
problems, and these problems are of the same level of difficulty.
Generalized projection operators
are constructive and they obey computer processing at least in spaces
$l^p, L^p $  and Sobolev spaces $ W^p_m$  where  $\infty > p >1.$
This means that generalized projection operators in Banach spaces can be
viewed
as a natural generalization of the metric projection operators in
Hilbert spaces .\\
\\
Results of this paper were  presented at the 1993 SIAM annual meeting
in Philadelphia, July  12-16.


\begin{thebibliography}{99}

\bibitem{al4}  Ya.I. Alber, The regularization method for variational
inequalities with  nonsmooth
unbounded  operator in Banach space, {\it Appl. Math. Letters,\/}
6 (1993), 63-68.

\bibitem{al1} Ya.I. Alber, Recurrence relations and variational
inequalities,
{\it Soviet Math. Dokl.\/,}  27 (1983), 511-517.

\bibitem{al0} Ya.I. Alber, On the solution of equations and variational
inequalities
 with maximal monotone operators, {\it Soviet Math. Dokl.\/,}  20 (1979),
 871-876.

\bibitem{al2} Ya.I. Alber, Global version of Bjornestal's estimate for
metric projection
operator in Banach space  (to appear).

\bibitem{al3}  Ya.I. Alber, On penalty method for variational
inequalities with  nonsmooth
unbounded  operator in Banach space  (to appear).

\bibitem{aln9} Ya.I. Alber and A.I. Notik, Parallelogram inequalities in
Banach
spaces and some properties of the duality mapping, {\it Ukranian Math.
J.\/},
40 (1988), 650-652.

\bibitem{aln1} Ya.I. Alber and A.I. Notik, Geometric properties of
Banach spaces and approximate
methods for solving nonlinear operator equations, {\it Soviet Math.
Dokl.\/,}
 29 (1984), 611-615 .

\bibitem{aln2} Ya.I. Alber and A.I. Notik, On some estimates for
projection operator in Banach
space (to appear).

\bibitem{aln3} Ya.I. Alber and A.I. Notik, On minimization of
functionals and
solution of variational inequalities in Banach spaces, {\it Soviet Math.
Dokl.,\/}
34 (1986), 296-300.

\bibitem{aln8} Ya.I. Alber and A.I. Notik, Geometry of Banach spaces and
approximation methods for solving nonlihear operator equations,
Preprint,
Instit. Radio Phisics Researches, Gorky, 1983.

\bibitem{aus}  A. Auslender,  {\it Optimization. Methodes numeriques.
Maitrisede Mathematiques et Applications Fundamentales  \/},  Masson,
Paris-New York-Barcelona, 1976.

\bibitem{bj}  B.O.Bjornestal, Local Lipschitz continuity of the metric
projection
operator,  {\it  Approximation Theory, \/} Stefan Banach
Internat. Math. Center Publication, Warsaw, 4 (1979), 43-53.

\bibitem{breg} L.M. Bregman, The method of successive projection  for
finding a common point
of convex sets, {\it Soviet Math. Dokl.,\/}  6 (1965), 688-692.

\bibitem{br} F.E. Browder, Nonlinear operators and nonlinear equations
in
Banach space, {\it Proc. Symp. Pure . Math. 18, Part II, Amer. Math.
Sos.,\/}
Providence, 1979.

\bibitem{bc} D. Butnariu and Y. Censor,  Strong convergence of almost
simultaneous
 block-iterative projection methods in Hilbert spaces,  (to appear).

\bibitem{cea}  J. Cea, {\it Optimization: Theorie et algorithmes.
Methodes
Mathematiques de l'Informatique, 2  \/,} Dunod, Paris, 1971.

\bibitem{ct} P.L. Combettes and H.J. Trussell, The method of successive
projection  for finding
a common point of sets in metric spaces,
{\it J. of Optimization Theory and Applications,\/}
67 (1990), 487-507.

\bibitem{d}  F. Deutsch,  The method of alternating orthogonal projections,
{\it Approximation Theory, Spline Functions and Applications, \/}
Kluwer Academic Publishers, S.P. Singh, ed., 105 - 121, 1992.

\bibitem{dis} J. Diestel',  The Geometry of Banach Spaces,
{\it Lecture Notes Math.,485, \/} Springer, 1975.

\bibitem{dunn} J.C.  Dunn,  Clobal and asymptotic rate estimates for a
class of
projected  gradient processes, {\it SIAM  J. Control  and
Optimization,\/}
19 (1981), 368-400.

\bibitem{fg}  T. Figiel,  On the moduli of convexity anm smoothness,
{\it Studia Mathematica, \/} 56 (1976), 121-155.

\bibitem{gr}  K. Goebel and S. Reich, {\it  Uniform convexity,
hyperbolic
geometry and nonexpansive mapping,\/} New York, N.Y., M.Dekker, 1984.

\bibitem{gpr} L.G. Gubin, B.T.Polyak and  E.V. Raik,
The method of  projections  for finding the
 common point of  convex sets, {\it USSR Computational Mathematics and
  Mathematical Physics,\/}  7 (1967), 1-24.

\bibitem{h} R.B. Holmes,
{\it A Course on Optimization and Best Approximation, \/} Springer-
Verlag,
Berlin - Heidelberg - New York, 1972.

\bibitem{jy} C.R. Joy and J.C. Yao, Algorithm for generealized
multivalued variational inequalities in Hilbert spaces, {\it Computer
Math. Applic.,\/} 25 (1993), 7-16.

\bibitem{ls} T.C. Lim and R. Smarzewski, On best approximation and
coapproximation in $L^p$ spaces, {\it Progress in Approximation
Theory\/},
Academic Press, Inc. Harcourt Brace Jovanovich Publishers, P.Nevai and
A.Pinkus, ed., 625 - 628, 1991.

\bibitem{lhh} T.C. Lim, H.K. Xu and Z.B. Xu, An $L^p$ inequality and its
applications to fixed point theory and approximation theory,
{\it Progress in Approximation Theory\/}, Academic Press, Inc. Harcourt
Brace Jovanovich Publishers, P.Nevai and A.Pinkus, ed., 609-624, 1991.

\bibitem{lz} J. Lindenstrauss and L. Tzafriri, {\it Classical Banach Spaces II
, \/} Springer-Verlag, Berlin-Heidelberg-New York, 1979.

\bibitem{li} J.-L. Lions,
{\it Quelques methodes de resolution des problemes aux limites non
lineaires,\/}  Dunod Gauthier-Villars, Paris, 1969.

\bibitem{ml} V.D. Mil'man,
Geometric theory of Banach Spaces. Part II. Geometry of the unit sphere,
{\it Russian Math. Surveys, \/} 26 (1971), No.6, 79-163.

\bibitem{n} M.Z. Nashed, A decomposition relative to convex sets,
{\it Proc. Am. Math. Soc.,\/} 19 (1968), 782-786.

\bibitem{no} M.A. Noor, Quasi variational inequalities,  {\it Appl. Math.
Lett.,\/} 1 (1988), 367-370.

\bibitem{not} A.I. Notik,  Properties of a duality mapping with a
scale function, {\it Soviet Math.,\/}
29 (1985), 96-98.

\bibitem{pet}  W.V. Petryshyn,  Projection methods in nonlinear
numarical
 functional analysis, {\it  J. Math. Mech. ,\/} 17 (1967), 353-372.

\bibitem{sm} B. Prus  and R. Smarzewski, Strongly unique best
approxination
 and centers in uniformly convex spaces,  {\it J.Math. Anal. Appl.,\/}
121 (1987), 10-21.

\bibitem{r} S. Reich, Extension problems for accretive sets in Banach
spaces,
 {\it J. Functional Anal.,\/} 26 (1977), 378-395.

\bibitem{shi} P. Shi, Equivalence of variational inequalities with
Wiener-Hopf equations, {\it Proc. Amer. Math. Soc.,\/} 111 (1991), 339-346.

\bibitem{vnc} A.A. Vladimirov, Yu.E. Nesterov and Yu.N. Chekanov, Uniformly
convex functionals, {\it Vestnik Moscov. Univ. Ser.15, Vychisl. Mat.
Kibernet.,\/} 1978, No.3, 12-23.

\bibitem{vn} J. Von-Neumann, On rings of operators. Reduction theory,
{\it Annals of Mathematics,\/} 50 (1949), 401-485.

\bibitem{Yu} V.V. Yurgelas, Some geometric characteristics of Banach
spaces and
accretive operators, {\it Izv.Visch. Uchebn. Zaved. Math.,\/}
5 (1982), 63-69.

\bibitem{zr} Z.B. Xu and G.F. Roach, Charateristic inequalities of
uniformly convex
and uniformly smooth Banach Spaces, {\it J. Math. Anal. Appl.,\/}
157 (1991), 189-210.

\bibitem{z} E.H.Zarantonello, Projections on convex sets in Hilbert space
and spectral theory, {\it Contributions to Nonlinear Functional Analysis,\/}
Academic Press, New York - London, E.H.Zarantonello, ed., 1972.

\bibitem{zd} E. Zeidler, {\it Nonlinear Functional Analysis and its
Applications 11/B: Nonlinear
Monotone Operators, \/} Springer-Verlag, New-York Inc., 1990.

\end{thebibliography}
\end{document}